\documentclass{amsart}

\usepackage{amsmath,amscd,amssymb,amsthm}
\usepackage{enumerate,mathrsfs}
\usepackage{geometry}
\usepackage{hyperref}
\usepackage{url}
\usepackage{ifthen}
\usepackage[all]{xy} \xyoption{2cell} \UseAllTwocells
\usepackage[mathscr]{euscript}

\numberwithin{equation}{section}

\theoremstyle{plain}
\newtheorem{thm_}[equation]{Theorem}
\newtheorem{lemma_}[equation]{Lemma}
\newtheorem{prop_}[equation]{Proposition}
\newtheorem{cor_}[equation]{Corollary}
\newtheorem{eg_}[equation]{Example}
\newtheorem{con_}[equation]{Conjecture}
\newtheorem*{cons_}{Conjecture}

\theoremstyle{definition}
\newtheorem{thmu_}[equation]{Theorem}
\newtheorem*{thmus_}{Theorem}
\newtheorem{propu_}[equation]{Proposition}
\newtheorem*{propus_}{Proposition}
\newtheorem{coru_}[equation]{Corollary}
\newtheorem*{corus_}{Corollary}
\newtheorem{lemu_}[equation]{Lemma}
\newtheorem*{lemus_}{Lemma}
\newtheorem{egu_}[equation]{Example}
\newtheorem*{egus_}{Example}
\newtheorem{def_}[equation]{Definition}
\newtheorem*{defs_}{Definition}
\newtheorem{rk_}[equation]{Remark}
\newtheorem*{rks_}{Remark}
\newtheorem{ex_}[equation]{Remark}

\newcommand{\prop}[1]{\begin{prop_}#1\end{prop_}}

\newcommand{\rk}[1]{\begin{rk_}#1\end{rk_}}

\newcommand{\cor}[1]{\begin{cor_}#1\end{cor_}}

\newcommand{\pf}[1]{\begin{proof}#1\end{proof}}

\DeclareMathOperator{\Gal}{Gal}
\DeclareMathOperator{\Hom}{Hom}

\DeclareMathOperator{\Ext}{Ext}

\DeclareMathOperator{\Spec}{Spec}

\DeclareMathOperator{\im}{im}
\DeclareMathOperator{\Ob}{Ob}

\DeclareMathOperator{\Br}{Br}

\DeclareMathOperator{\Pic}{Pic}

\DeclareMathOperator{\KD}{KD}

\newcommand{\HH}{\mathbb H}

\newcommand{\RR}{\mathbb R}
\newcommand{\ZZ}{\mathbb Z}
\newcommand{\bfD}{\mathbf D}%
\newcommand{\bfG}{\mathbf G}%
\newcommand{\bfR}{\mathbf R}

\newcommand{\cA}{\mathcal A}%
\newcommand{\cB}{\mathcal B}%
\newcommand{\cC}{\mathcal C}%

\newcommand{\La}{\Lambda}
\newcommand{\tm}{\times}%

\newcommand{\ol}{\overline}

\newcommand{\ra}{\rightarrow}

\newcommand{\Ra}{\Rightarrow}

\newcommand{\xra}{\xrightarrow}

\newcommand{\is}[2]{\xymatrix@-4mm{#1 \ar[r]^-{\sim} & #2 }}
\newcommand{\mis}[2]{\xymatrix@-2mm{#1 \ar[r]^-{\sim} & #2 }}

\newcommand{\dra}[4]{\xymatrix@-4mm{#1 \ar@<.5ex>[r]^-{#3} \ar@<-.5ex>[r]_-{#4}& #2 }}
\newcommand{\era}[5]{\xymatrix@-4mm{#1 \ar[r] &#2 \ar@<.5ex>[r]^-{#4} \ar@<-.5ex>[r]_-{#5}& #3 }}

\newcommand{\tu}[1]{\text{\upshape #1}}

\newcommand{\AAb}{{Ab}}
\newcommand{\MMod}{{Mod}}

\newcommand{\eq}[1]{\begin{equation}#1\end{equation}}
\newcommand{\eqn}[1]{\begin{equation*}#1\end{equation*}}

\newcommand{\gan}[1]{\begin{gather*}#1\end{gather*}}
\newcommand{\al}[1]{\begin{align}#1\end{align}}
\newcommand{\aln}[1]{\begin{align*}#1\end{align*}}

\newcommand{\enmt}[1]{\begin{enumerate}#1\end{enumerate}}

\newcommand{\aci}[1]{\ar@{^(->}[#1]|-{/}}
\newcommand{\coaci}[1]{\ar@{_(->}[#1]|-{/}}
\newcommand{\aoi}[1]{\ar@{^(->}[#1]|-{\circ}}
\newcommand{\coaoi}[1]{\ar@{_(->}[#1]|-{\circ}}
\makeatletter
\def\citet@url@sp{https://stacks.math.columbia.edu/}
\def\citet@bib@sp{stacks-project}
\def\citet@url@kd{https://kerodon.net/}
\def\citet@bib@kd{kerodon}
\newcommand{\citet@tag}[2]{\href{#2tag/#1}{#1}}
\newcommand{\citet@taglist}[2]{%
 \def\@citet@e{}%
 \def\@citet@tag@n{0}
 \@for\@citet@tag:=#1\do{%
  \edef\@citet@tag@n{\the\numexpr\@citet@tag@n + 1}%
 }%
 \def\@citet@tags{%
  \def\@citet@tag@i{0}%
  \@for\@citet@tag:=#1\do{%
   \edef\@citet@tag@i{\the\numexpr\@citet@tag@i + 1}%
   \ifthenelse{\@citet@tag@i > 1}{
    \ifthenelse{\@citet@tag@i = \@citet@tag@n}{
     \citet@seplast%
    }{%
     \citet@sep%
    }%
   }{}%
   \citet@entry{\@citet@tag}{#2}
  }%
 }%
 \ifthenelse{\@citet@tag@n > 1}{
  \def\@citet@Tag{Tags}%
 }{%
  \def\@citet@Tag{Tag}%
 }%
 \@citet@Tag~\@citet@tags
}
\newcommand{\citet@sep}{, }
\newcommand{\citet@seplast}{ and }
\newcommand{\citet@entry}[2]{\citet@tag{#1}{#2}}
\let\@old@cite\cite
\renewcommand{\cite}[2][]{%
 \def\@citet@detail{\citet@taglist{#1}{\@citet@url}}%
 \ifthenelse{\equal{#2}{sp}}{%
  \def\@citet@url{\citet@url@sp}%
  \def\@citet@bib{\citet@bib@sp}%
 }{\ifthenelse{\equal{#2}{kd}}{%
  \def\@citet@url{\citet@url@kd}%
  \def\@citet@bib{\citet@bib@kd}%
 }{
  \def\@citet@detail{#1}%
  \def\@citet@bib{#2}%
 }}%
 \ifthenelse{\equal{#1}{}}{%
  \@old@cite{\@citet@bib}%
 }{%
  \@old@cite[\@citet@detail]{\@citet@bib}%
 }%
}
\makeatother

\newcommand{\etale}{{\'etale}}

\newcommand{\Grot}{{Grothendieck}}

\newcommand{\desc}{\tu{desc}}

\newcommand{\hdesc}[2]{{#1\tu{-}\desc{\ifthenelse{\equal{#2}{}}{}{_#2}}}}

\begin{document}
\title[Commutative diagrams of low term exact sequences of Ext]{
 On  commutative diagrams consisting of low term exact sequences}
\author[C. Lv]{Chang Lv}
\address{State Key Laboratory of Information Security\\
Institute of Information Engineering\\
Chinese Academy of Sciences\\
Beijing 100093, P.R. China}
\email{lvchang@amss.ac.cn}
\date{\today}
\thanks{This work was supported by
 National Natural Science Foundation of China (Grant No. 11701552).
}
\begin{abstract}
We establish several useful commutative diagrams
 consisting of low term exact sequences
 attached to {\Grot} spectral sequences,
 which extends and
 integrates the previous ones appeared in literature
 such as Alexei~N. Skorobogatov [Beyond the {M}anin obstruction, Invent. Math.
  (1999)], and [On the elementary obstruction to the existence of rational
  points, Mathematical Notes (2007)].
Parts of the diagrams was frequently
 used in local-global principle to rational points.
\end{abstract}
\maketitle

\section{The commutative diagram} \label{sec_main}
Suppose that $\Phi: \cA\ra\cB$ and $\Psi_t: \cB\ra\cC$ are left exact additive
 functors between abelian categories, $t=1,2,3$.
Assume that  $\cA$ and $\cB$ have enough injectives and $\Psi_t$ takes injectives to
 $\Phi$-acyclics.
Then for any $A\in\Ob(\cA)$, we have the {\Grot} spectral sequence
\eq{
\tag*{$(S)_{\Phi, \Psi_t, A}$} \label{eq_groth}
^tE_2^{p,q}=(R^p\Psi)(R^q\Phi)A\Ra \ ^tE^{p+q}=R^{p+q}(\Psi\Phi)A
}
 and the low term exact sequence
\eq{
\tag*{$(E)_{\Phi, \Psi_t, A}$} \label{eq_low}
0\ra\ ^tE_2^{1,0}\ra\ ^tE^1\ra\ ^tE_2^{0,1}\ra\ ^tE_2^{2,0}\ra\ ^tE_1^2\ra\
 ^tE_2^{1,1}\ra\ ^tE_2^{3,0}
}
 attached to them, where $^tE_1^2=\ker(^tE^2\ra\ ^tE_2^{0,2})$.
Let  $\bfD^+(\cA)$, $\bfD^+(\cB)$, $\bfD^+(\cC)$ be the  corresponding
 derived category of complexes
 bounded below, $\bfR\Psi$, $\bfR\Psi_t$  the corresponding derived functor and
  $\RR^\bullet$ (or  $\HH^\bullet$)  the hyppercohomology functor.
\prop{\label{prop_e2terms}
With the previous notation, suppose that there are morphism of functors
 $u: \Psi_1\ra \Psi_2$ and $v: \Psi_2\ra \Psi_3$ such that for any $F\in\bfD^+(\cB)$,
\eq{\label{eq_psi_tri}
\xymatrix{
\bfR\Psi_1(F)\ar[r]^-{\bfR u(F)} &\bfR\Psi_2(F)\ar[r]^-{\bfR v(F)}
 &\bfR\Psi_3(F)\ar[r] &\bfR\Psi_1(F)[1]
}}
 is a distinguished triangle functorial in $F$.

\enmt{[\upshape (i)]
\item \label{it_long}
We have the long exact sequence
\eqn{
\dots\ra\ ^tE_2^{i,0}\ra\ ^tE_{\le1}^i\ra\ ^tE_2^{i-1,1}\ra\ ^tE_2^{i+1,0}\ra\dots
}
 where $^tE_{\le1}^i=\RR^i\Psi_t(\tau_{\le1}\bfR\Phi(A))$.
\item \label{it_com}
We have the following commutative diagram with exact rows and columns
\eq{\label{eq_e2terms}
\xymatrix{
^1E_2^{1,0}\ar[r]\ar[d] &^1E^1\ar[r]\ar[d] &^1E_2^{0,1}\ar[r]\ar[d]
 &^1E_2^{2,0}\ar[r]\ar[d] &^1E_1^2\ar[d] \\
^2E_2^{1,0}\ar[r]\ar[d] &^2E^1\ar[r]\ar[d] &^2E_2^{0,1}\ar[r]\ar[d]
 &^2E_2^{2,0}\ar[r]\ar[d] &^2E_1^2\ar[d] \\
^3E_2^{1,0}\ar[r]\ar[d] &^3E^1\ar[r]\ar[d] &^3E_2^{0,1}\ar[r]\ar[d]
 &^3E_2^{2,0}\ar[r]\ar[d] &^3E_1^2\ar[d] \\
^1E_2^{2,0}\ar[r] &^1E_1^2\ar[r] &^1E_2^{1,1}\ar[r] &^1E_2^{3,0}\ar[r]
 &^1E_{\le1}^3
}}
 where $^1E_{\le1}^3$ fits into the exact sequence
\eq{\label{eq_e3}
0\ra\ ^1E_1^2\ra\ ^1E^2\ra\ ^1E_2^{0,2}\ra\ ^1E_{\le1}^3\ra\ ^1E^3,
}
 the rows are parts of the low term exact sequences \ref{eq_low} attached to
 the spectral sequences \ref{eq_groth} with $t$ numbered on the left upper corner
 of each object, and
 the columns are induced by taking cohomology at $1$ of \eqref{eq_psi_tri} in which
 $F$ is substituted with $\tau_{[0]}\bfR\Phi(A)$, $\tau_{\le1}\bfR\Phi(A)$,
 $\tau_{[1]}\bfR\Phi(A)$, $\tau_{[0]}\bfR\Phi(A)[1]$, $\tau_{\le1}\bfR\Phi(A)[1]$,
 respectively.

\item \label{it_chase}
 Suppose we are given $\beta\in\ ^3E^1$ and
 $\gamma\in\ ^2E_2^{0,1}$ such that they map to the same element in $^3E_2^{0,1}$.
 Then there exists $\alpha\in\ ^1E_2^{2,0}$ such that $\alpha,\ \beta$ map to
 the same element in $^1E_1^2$ and $-\alpha,\ \gamma$ map to the same element
 in $^2E_2^{2,0}$. In other words, we have the zig-zag diagram
\eqn{\xymatrix@-6mm{
 & & & & & -\alpha\ar[d] &^1E_2^{2,0}\ar[dd] \\
 & & &\gamma\ar[rr]\ar[dd] & &c & \\
 & & & &^2E_2^{0,1}\ar[rr]\ar[dd] & &^2E_2^{2,0} \\
 &\beta\ar[rr]\ar[dd] & &a & & & \\
 & &^3E^1\ar[rr]\ar[dd] & &^3E_2^{0,1} & & \\
 \alpha\ar[r] &b & & & & & \\
 ^1E_2^{2,0}\ar[rr] & &^1E_1^2 & & & &
}}

\item \label{it_chase2}
 The statement of \eqref{it_chase} is also correct if we move our focus one step
 right. That is, we are given $\beta\in\ ^3E_2^{0,1}$ and $\gamma\in\ ^2E_2^{2,0}$,
 and so on.
}}
\pf{
The proof widely extends  \cite[Lem. 3]{skorobogatov99beyond}
 and \cite[Prop. 1.1]{skorobogatov2007eltob}.
For  any $F\in\bfD^+(\cA)$, the truncation functors determine the distinguished triangle
\eq{\label{eq_trun}
\tau_{\le0}F\ra F\ra \tau_{\ge1}F\ra (\tau_{\le0}F)[1].
}
Note that $\bfR\Psi_t$, $t=1,2,3$ are triangulated.
Along with the functorial distinguished triangles \eqref{eq_psi_tri}
 in which $F$ is substituted with $\tau_{\le0}F$, $F$ and $\tau_{\ge1}F$ respectively,
 we obtain the following commutative diagram
\eq{\label{eq_tri}
\xymatrix{
\bfR\Psi_1(\tau_{\le0}F)\ar[r]\ar[d] &\bfR\Psi_1(F)\ar[r]\ar[d]
 &\bfR\Psi_1(\tau_{\ge1}F)\ar[r]\ar[d] &\bfR\Psi_1(\tau_{\le0}F)[1]\ar[r]\ar[d]
 &\bfR\Psi_1(F)[1]\ar[d] \\
\bfR\Psi_2(\tau_{\le0}F)\ar[r]\ar[d] &\bfR\Psi_2(F)\ar[r]\ar[d]
 &\bfR\Psi_2(\tau_{\ge1}F)\ar[r]\ar[d] &\bfR\Psi_2(\tau_{\le0}F)[1]\ar[r]\ar[d]
 &\bfR\Psi_2(F)[1]\ar[d] \\
\bfR\Psi_3(\tau_{\le0}F)\ar[r]\ar[d] &\bfR\Psi_3(F)\ar[r]\ar[d]
 &\bfR\Psi_3(\tau_{\ge1}F)\ar[r]\ar[d]
 &\bfR\Psi_3(\tau_{\le0}F)[1]\ar[r]\ar[d] &\bfR\Psi_3(F)[1]\ar[d] \\
\bfR\Psi_1(\tau_{\le0}F)[1]\ar[r] &\bfR\Psi_1(F)[1]\ar[r] &\bfR\Psi_1(\tau_{\ge1}F)[1]
 \ar[r] &\bfR\Psi_1(\tau_{\le0}F)[2]\ar[r] &\bfR\Psi_1(F)[2]
}}
Let $F=\tau_{\le1}\bfR\Phi(A)=\tau_{[0,1]}\bfR\Phi(A)$ and taking cohomology at $1$,
 the diagram becomes
\eq{\label{eq_rr1}
\scriptsize\xymatrix@-5mm{
\RR^1\Psi_1(\tau_{[0]}\bfR\Phi(A))\ar[r]\ar[d] &\RR^1\Psi_1(\tau_{\le1}\bfR\Phi(A))
 \ar[r]\ar[d] &\RR^1\Psi_1(\tau_{[1]}\bfR\Phi(A))\ar[r]\ar[d]
 &\RR^1\Psi_1(\tau_{[0]}\bfR\Phi(A))[1]\ar[r]\ar[d]
 &\RR^1\Psi_1(\tau_{\le1}\bfR\Phi(A))[1]\ar[d] \\
\RR^1\Psi_2(\tau_{[0]}\bfR\Phi(A))\ar[r]\ar[d] &\RR^1\Psi_2(\tau_{\le1}\bfR\Phi(A))
 \ar[r]\ar[d] &\RR^1\Psi_2(\tau_{[1]}\bfR\Phi(A))\ar[r]\ar[d]
 &\RR^1\Psi_2(\tau_{[0]}\bfR\Phi(A))[1]\ar[r]\ar[d]
 &\RR^1\Psi_2(\tau_{\le1}\bfR\Phi(A))[1]\ar[d] \\
\RR^1\Psi_3(\tau_{[0]}\bfR\Phi(A))\ar[r]\ar[d] &\RR^1\Psi_3(\tau_{\le1}\bfR\Phi(A))
 \ar[r]\ar[d] &\RR^1\Psi_3(\tau_{[1]}\bfR\Phi(A))\ar[r]\ar[d]\ar[r]
 &\RR^1\Psi_3(\tau_{[0]}\bfR\Phi(A))[1]\ar[r]\ar[d]
 &\RR^1\Psi_3(\tau_{\le1}\bfR\Phi(A))[1]\ar[d] \\
\RR^1\Psi_1(\tau_{[0]}\bfR\Phi(A))[1]\ar[r] &\RR^1\Psi_1(\tau_{\le1}\bfR\Phi(A))[1]
 \ar[r] &\RR^1\Psi_1(\tau_{[1]}\bfR\Phi(A))[1]\ar[r]
 &\RR^1\Psi_1(\tau_{[0]}\bfR\Phi(A))[2]\ar[r] &\RR^1\Psi_1(\tau_{\le1}\bfR\Phi(A))[2]
}}
 with exact rows and columns.

We now identify the objects appearing in \eqref{eq_rr1} with the ones in
 \eqref{eq_e2terms}.
Clearly for any $t=1,2,3$ and $i, j, k\in\ZZ$ with $j, i-j+k\ge0$, we have
\eqn{
\RR^i\Psi_t(\tau_{[j]}\bfR\Phi(A))[k] = (R^{i-j+k}\Psi_t)(R^j\Phi)A =\ ^tE_2^{i-j+k, j}.
}
It remains to identify $\RR^1\Psi_t(\tau_{\le1}\bfR\Phi(A))$  with $^tE^1$ and
 $\RR^1\Psi_t(\tau_{\le1}\bfR\Phi(A))[1]$ with $^tE_1^2$.

For any $F$ consider the distinguished triangle
\eqn{
\tau_{\le1}F\ra F\ra \tau_{\ge2}F\ra (\tau_{\le1}F)[1]
}
 and note that $\Psi_t(\tau_{\ge2}F)$ is acyclic in $0$ and $1$.
Then we have the long exact sequences
\al{\label{eq_long}
\RR^0\Psi_t(\tau_{\ge2}F)\ra &\RR^1\Psi_t(\tau_{\le1}F)\ra \RR^1\Psi_t(F)\ra
 \RR^1\Psi_t(\tau_{\ge2}F) \nonumber\\
 \ra &\RR^2\Psi_t(\tau_{\le1}F)\ra \RR^2\Psi_t(F)\ra
 \RR^2\Psi_t(\tau_{\ge2}F) \nonumber\\
 \ra &\RR^3\Psi_t(\tau_{\le1}F)\ra \RR^3\Psi_t(F)
}
 where $\RR^i\Psi_t(\tau_{\ge2}F)=0$ with $t=1,2,3$ and $i=0,1$.
Now we take $F=\bfR\Phi(A)$. It follows that
\eqn{
\RR^1\Psi_t(\tau_{\le1}\bfR\Phi(A)) = \RR^1\Psi_t(\bfR\Phi(A)) =\ ^tE^1
}
 where the last equality follows from the isomorphism of functors
\eq{\label{eq_iso_der}
(\bfR\Psi_t)(\bfR\Phi)\cong \bfR(\Psi_t\Phi).
}

In a same manner,  \eqref{eq_long} and \eqref{eq_iso_der} yield
\aln{
\RR^1\Psi_t(\tau_{\le1}\bfR\Phi(A))[1] &= \RR^2\Psi_t(\tau_{\le1}\bfR\Phi(A))\\
 &=\ker\left(\RR^2\Psi_t(\bfR\Phi(A))\ra \RR^2\Psi_t(\tau_{\ge2}\bfR\Phi(A))\right)\\
 &=\ker\left(\RR^2(\Psi_t\Phi)A\ra \RR^2\Psi_t(\tau_{\ge2}\bfR\Phi(A))\right).
}
Then  the low term exact sequence  attached to the hyppercohomology spectral sequence
 \cite[Appendix C (g)]{milne80etale}
\eqn{
E_2^{p,q}=R^pg(H^q(F))\Ra \RR^{p+q}g(F)
}
gives the isomorphism when $F=\tau_{\ge2}\bfR\Phi(A)$ and $g=\Psi_t$
\eqn{
\RR^2\Psi_t(\tau_{\ge2}\bfR\Phi(A))\cong \Psi_t(R^2\Phi(A)),
}
 since in this case $E_2^{p,q}=0$ for all $p$ and $q=0,1$.
Thus\eqn{
\RR^1\Psi_t(\tau_{\le1}\bfR\Phi(A))[1] =
 \ker\left(\RR^2(\Psi_t\Phi)A\ra \Psi_t(R^2\Phi(A))\right) =\ ^tE_1^2.
}

Similarly,
\eqn{
\RR^1\Psi_1(\tau_{\le1}\bfR\Phi(A))[2]=\RR^3\Psi_1(\tau_{\le1}\bfR\Phi(A))=
 \ ^1E_{\le1}^3
}
and the exact sequence \eqref{eq_e3} is also  deduced from
\eqref{eq_long} and \eqref{eq_iso_der}.
This completes the identification of objects.

Finally, in  diagram \eqref{eq_rr1}, the identification of the vertical arrows
 is clear.
For that of the horizontal ones, it follows  from a general fact for such spectral
 sequences. See, for example, \cite[Appendix B]{skorobogatov99beyond}, which shows
 that $\RR^1\Psi_t(\tau_{\le1}\bfR\Phi(A))\ra \RR^1\Psi_t(\tau_{[1]}\bfR\Phi(A))$
 is exactly the edge map $^tE_1\ra\ ^tE_2^{0,1}$.
This completes the proof of \eqref{it_com} as well as \eqref{it_long}.

Consider the subdiagram of \eqref{eq_tri}
\eq{\label{eq_sub}
\xymatrix{
\bfR\Psi_2(F)\ar[r]\ar[d] &\bfR\Psi_2(\tau_{\ge1}F)\ar[r]\ar[d]
 &\bfR\Psi_2(\tau_{\le0}F)[1]\ar[r]\ar[d] & \\
\bfR\Psi_3(F)\ar[r]\ar[d] &\bfR\Psi_3(\tau_{\ge1}F)\ar[r]\ar[d]
 &\bfR\Psi_3(\tau_{\le0}F)[1]\ar[r]\ar[d] & \\
\bfR\Psi_1(F)[1]\ar[r]\ar[d] &\bfR\Psi_1(\tau_{\ge1}F)[1]
 \ar[r]\ar[d] &\bfR\Psi_1(\tau_{\le0}F)[2]\ar[r]\ar[d] & \\
 & & &
}}
 whose rows and columns are all distinguished triangles.
Since up to an isomorphism, every distinguished triangle in a derived category
 arises from some  short exact sequence of complexes
 \cite[Chap. IV.2 8. Prop.]{gelfand2013methods},
 we may view \eqref{eq_sub} as a commutative diagram consisting of three rows and
 three columns of short exact sequences of complexes.
Then the result follows from  \cite[Lem. 4.3.2]{torsor} by taking $i=1$.
This completes the proof of \eqref{it_chase}.

\eqref{it_chase2} is similar as  \eqref{it_chase}. The proof is complete.
}

Next we describe a variant of Proposition \ref{prop_e2terms}.
\prop{\label{prop_var_cone}
Keeping assumptions in {\upshape Proposition \ref{prop_e2terms}}, suppose that
 there are $A\in\bfD^+(\cA)$ and $B\in\bfD^+(\cB)$ with
 a morphism \eqn{f: B\ra \tau_{\le1}\bfR\Phi(A).}
Let $\Delta=\Delta(\Phi, A, B, f)$ be the cone of $-f[1]$.
Denote $^tF^p=\RR^p\Psi_tB$ and $^tG^p=\RR^p\Psi_t\Delta$.
Then we have the long exact sequence
\eqn{
\dots\ra F^i\ra\ ^tE_{\le1}^i\ra G^{i-1}\ra F^{i+1}\ra\dots
}
 where $^tE_{\le1}^i=\RR^i\Psi_t(\tau_{\le1}\bfR\Phi(A))$
 and the following commutative diagram with exact rows and columns
\eq{\label{eq_var_cone}
\xymatrix{
^1F^1\ar[r]\ar[d] &^1E^1\ar[r]\ar[d] &^1G^0\ar[r]\ar[d]
 &^1F^2\ar[r]\ar[d] &^1E_1^2\ar[d] \\
^2F^1\ar[r]\ar[d] &^2E^1\ar[r]\ar[d] &^2G^0\ar[r]\ar[d]
 &^2F^2\ar[r]\ar[d] &^2E_1^2\ar[d] \\
^3F^1\ar[r]\ar[d] &^3E^1\ar[r]\ar[d] &^3G^0\ar[r]\ar[d]
 &^3F^2\ar[r]\ar[d] &^3E_1^2\ar[d] \\
^1F^2\ar[r] &^1E_1^2\ar[r] &^1G^1\ar[r] &^1F^3\ar[r] & ^1E_{\le1}^3
}}
 where $^1E_{\le1}^3$ fits into the exact sequence
\eqn{
0\ra\ ^1E_1^2\ra\ ^1E^2\ra\ ^1E_2^{0,2}\ra\ ^1E_{\le1}^3\ra\ ^1E^3.
}

Moreover, similar statements as {\upshape \eqref{it_chase}, \eqref{it_chase2}} in
 {\upshape Proposition \ref{prop_e2terms}}  hold.
That is, if we are given $\beta\in\ ^3E^1$
 $($resp. $^3G^0)$ and $\gamma\in\ ^2G^0$ $($resp. $^2F^2)$, then the corresponding
 zig-zag
 diagram as in   {\upshape Proposition \ref{prop_e2terms}}  is correct.
}
\pf{
The same as Proposition \ref{prop_e2terms}  except that in the diagram
 \eqref{eq_rr1} we replace the distinguished triangle
\eq{\label{eq_tri_basic}
\tau_{[0]}\bfR\Phi(A)\ra \tau_{\le1}\bfR\Phi(A)\ra \tau_{[1]}\bfR\Phi(A)\ra
 \tau_{[0]}\bfR\Phi(A)[1]
}
by \eqn{
B\ra \tau_{\le1}\bfR\Phi(A)\ra \Delta[-1]\ra B[1].
}
}

\rk{
\enmt{[(a)]
\item Obviously, Proposition \ref{prop_e2terms} is the special case of
 Proposition \ref{prop_var_cone} where $f$ is the canonical map
 $\tau_{[0]}\bfR\Phi(A)\ra \tau_{\le1}\bfR\Phi(A)$.
\item If we replace the distinguished triangle \eqref{eq_tri_basic}
by \eqn{
\tau_{[0]}\bfR\Phi(A)\ra \bfR\Phi(A)\ra \tau_{\ge1}\bfR\Phi(A)\ra
 \tau_{[0]}\bfR\Phi(A)[1]
}
we also have  the long exact sequence
\eqn{
\dots\ra\ ^tE_2^{i,0}\ra\ ^tE^i\ra\ ^tE_{\ge1}^i\ra\ ^tE_2^{i+1,0}\ra\dots
}
 where $^tE_{\ge1}^i=\RR^i\Psi_t(\tau_{\ge1}\bfR\Phi(A))$.
}}

\section{Applications}
Suppose that $f_*: \cA\ra\cB$ is a left exact additive functor between two abelian
 categories which has a left adjoint $f^*$.
Assume that  $\cA$ and $\cB$ has enough injectives and $f^*$ is exact.
For example, $f: X\ra Y$ is a morphism of topoi, and $\cA=\MMod(X,\La)$,
 $\cB=\MMod(Y, \La)$.
Let $M\in\Ob(\cB)$ and $N\in\Ob(\cA)$.
Then we have the {\Grot} spectral sequences
\eq{\label{eq_extspc}
^ME_2^{p,q}=\Ext_\cA^p(M, R^qf_*N)\Ra \Ext_\cB^{p+q}(f^*M,N).
}
For simplicity, we omit the category letter in $\Ext$'s if it does not cause a confusion.
\cor{\label{cor_extterms}
Let $C, A\in \Ob(\cB)$
 and $u\in\Ext^1(C, A)$ be the element representing the extension
\eq{\label{eq_abc}
0\ra A\xra{i} B\xra{j} C\ra 0.
}
Define \gan{
\Ext_1^2(f^*M, N)=\ker\left(\Ext^2(f^*M, N)\ra\Hom(M, R^2f_*N)\right), \\
\Ext_{\le1}^3(f^*M, N)=\Ext^3(M, \tau_{\le1}\bfR f_*N).
}

\enmt{[\upshape (i)]
\item \label{it_extcom}
We have the following commutative diagram with exact rows and columns
\eqn{
\xymatrix@-2mm{
\Ext^1(C, f_*N) \ar[r]\ar[d]_-{j^*}
 &\Ext^1(f^*C, N) \ar[r]\ar[d]_-{f^*(j^*)}
 &\Hom(C, R^1f_*N) \ar[r]^-{^Cd_2^{0,1}}\ar[d]_-{j^*} &\Ext^2(C, f_*N)
 \ar[d]_-{j^*}\ar[r] &\Ext_1^2(f^*C, N)\ar[d]_-{f^*(j^*)} \\
\Ext^1(B, f_*N) \ar[r]\ar[d]_-{i^*}
 &\Ext^1(f^*B, N) \ar[r]\ar[d]_-{f^*(i^*)}
 &\Hom(B, R^1f_*N) \ar[r]^-{^Bd_2^{0,1}}\ar[d]_-{i^*} &\Ext^2(B, f_*N)
 \ar[d]_-{i^*}\ar[r] &\Ext_1^2(f^*B, N)\ar[d]_-{f^*(i^*)} \\
\Ext^1(A, f_*N) \ar[r]\ar[d]_-{u\cup-}
 &\Ext^1(f^*A, N) \ar[r]\ar[d]_-{f^*(u)\cup-}
 &\Hom(A, R^1f_*N) \ar[r]^-{^Ad_2^{0,1}}\ar[d]_-{u\cup-}
 &\Ext^2(A, f_*N)\ar[d]_-{u\cup-}\ar[r] &\Ext_1^2(f^*A, N)
 \ar[d]_-{f^*(u)\cup-} \\
\Ext^2(C, f_*N) \ar[r] &\Ext_1^2(f^*C, N) \ar[r] &\Ext^1(C, R^1f_*N)
 \ar[r]^-{^Cd_2^{1,1}}
 &\Ext^3(C, f_*N)\ar[r] &\Ext_{\le1}^3(f^*C, N)
}}
 where the rows are parts of the low term exact sequences  attached to
 $^AE_2^{p,q}$ and $^CE_2^{p,q}$ defined in \eqref{eq_extspc}.

\item \label{it_extchase}
 The statement of {\upshape Proposition \ref{prop_e2terms} \eqref{it_chase}}
 $($resp. \eqref{it_chase2}$)$ is also correct if we put $\beta,\gamma$ in the
 corresponding positions. That is, we are given $\beta\in\Ext^1(f^*A, N)$
 $($resp. $\Hom(A, R^1f_*N))$ and $\gamma\in\Hom(A, R^1f_*N)$ $($resp. $\Ext^2(B, f_*N))$,
 and so on.
}
}
\pf{
We shall use Proposition \ref{prop_e2terms}. Let
\eqn{
A\ra B\ra C\ra A[1]
}
 be the distinguished triangle in  $\bfD^+(\cB)$  determined by \eqref{eq_abc}.
Take $\cC=\AAb$,
  $\Psi_1 = \Hom(C, -)$,  $\Psi_2 = \Hom(B, -)$, $\Psi_3 = \Hom(A, -)$ and
 $\Phi = f_*$, which clearly satisfy the assumptions in Proposition \ref{prop_e2terms}
 (c.f. \cite[Thm. 10.7.4]{weibel1994homological}). Then the result follows.
}

Let $k$ be a field with characteristic $0$ and $\Gamma=\Gal(\ol k/k)$ where $\ol k$ is
 a fixed algebraic closure of $k$.
Let $p:X\ra \Spec k$ be a $k$-variety and $\ol X=X\tm_k \ol k$.
In Corollary \ref{cor_extterms} take $\cB$ be the category of discrete $\Gamma$-modules,
 $\cA$ the category of {\etale} sheaves on $X$.
We write $\Ext_k$ for $\Ext_\cB$, $\Ext_X$ for $\Ext_\cA$ and
\eqn{
\Ext_1^2(p^*T, \bfG_m)=\ker\left(\Ext_X^2(p^*T, \bfG_m)\ra\Hom_k(T, \Br\ol X)\right)
 \quad\text{ for any $\Gamma$-module $T$. }
}
Note that $\Ext_k^1(\ZZ, -)=H^1(k, -)$.
\cor{\label{cor_brextchase}
With previous notation, let $u\in H^1(k, M)$ be the element representing the extension
\eq{\label{eq_msz}
0\ra M\xra{i} S\xra{j} \ZZ\ra 0.
}
Define $H_{\le1}^3(X, \bfG_m)=H^3(X, \tau_{\le1}\bfR p_*\bfG_m)$.
\enmt{[\upshape (i)]
\item \label{cor_brextchase.it_brext}
We have the following commutative diagram with exact rows and columns
\eq{\label{eq_brextchase}
\xymatrix@-2mm{
H^1(k, \ol k[X]^\tm)\ar[r]\ar[d]_-{j^*} &\Pic X\ar[r]\ar[d]_-{p^*(j^*)}
  &\Pic\ol X^\Gamma\ar[r]\ar[d]_-{j^*}
 &H^2(k, \ol k[X]^\tm)\ar[r]\ar[d]_-{j^*} &\Br_1X\ar[d]_-{p^*(j^*)} \\
\Ext_k^1(S, \ol k[X]^\tm) \ar[r]^-{p^*}\ar[d]_-{i^*}
 &\Ext_X^1(p^*S, \bfG_m) \ar[r]\ar[d]_-{p^*(i^*)}
 &\Hom_k(S, \Pic\ol X)\ar[r]^-\partial\ar[d]_-{i^*}
 &\Ext_k^2(S, \ol k[X]^\tm)\ar[r]^-{p^*}\ar[d]_-{i^*}
 &\Ext_1^2(p^*S, \bfG_m)\ar[d]_-{p^*(i^*)} \\
\Ext_k^1(M, \ol k[X]^\tm) \ar[r]^-{p^*}\ar[d]_-{u\cup-}
 &\Ext_X^1(p^*M, \bfG_m) \ar[r]\ar[d]_-{p^*(u)\cup-}
 &\Hom_k(M, \Pic\ol X)\ar[r]^-\partial\ar[d]_-{u\cup-}
 &\Ext_k^2(M, \ol k[X]^\tm)\ar[r]^-{p^*}\ar[d]_-{u\cup-}
 &\Ext_1^2(p^*M, \bfG_m) \ar[d]_-{p^*(u)\cup-} \\
H^2(k, \ol k[X]^\tm)\ar[r] &\Br_1X\ar[r]^-r &H^1(k, \Pic\ol X)\ar[r]^-d
 &H^3(k, \ol k[X]^\tm)\ar[r] &H_{\le1}^3(X, \bfG_m)
}}
 where the top, middle and bottom row is a part of the low term exact sequences
 attached to the $E_2$ spectral sequences
 \eqn{\Ext_k^p(M, R^qp_*\bfG_m)\Ra \Ext_X^{p+q}(p^*M,\bfG_m),}
 \eqn{\Ext_k^p(S, R^qp_*\bfG_m)\Ra \Ext_X^{p+q}(p^*S,\bfG_m)} and
 \eqn{H^p(k, H^q(\ol X,\bfG_m))\Ra H^{p+q}(X,\bfG_m)} respectively.

\item \label{cor_brextchase.it_chase}
Let $\beta\in\Hom_k(M, \Pic\ol X)$ be such that  $u\cup\beta\in\im r$.
Then there exists $\alpha\in\Br_1X$ and $\gamma\in\Ext_k^2(S, \ol k[X]^\tm)$ such that
 $r(\alpha)=u\cup\beta$, $i^*(\gamma)=\partial(\beta)$ and
 $p^*(j^*)(-\alpha)=p^*(\gamma)$.
In other words, we have the diagram
\eqn{\xymatrix@-6mm{
 & & & & & -\alpha\ar[d] &\Br_1X\ar[dd] \\
 & & &\gamma\ar[rr]\ar[dd] & &p^*(\gamma) & \\
 & & & &\Ext_k^2(S, \ol k[X]^\tm)\ar[rr]\ar[dd] & &\Ext_1^2(p^*S, \bfG_m) \\
 &\beta\ar[rr]\ar[dd] & &\partial(\beta) & & & \\
 & &\Hom_k(M, \Pic\ol X)\ar[rr]\ar[dd] & &\Ext_k^2(M, \ol k[X]^\tm) & & \\
 \alpha\ar[r] &r(\alpha) & & & & & \\
 \Br_1X\ar[rr] & &H^1(k,\Pic\ol X) & & & &
}}}}
\pf{
Use Corollary \ref{cor_extterms}.
Take $f_*=p_*$, $N=\bfG_m$ and \eqref{eq_abc} to be \eqref{eq_msz}.
Then \eqref{cor_brextchase.it_brext} follows from the facts that for $p,q\ge0$,
 $R^qp_*\bfG_m=H^q(\ol X, \bfG_m)$
 and $H^p(X, \bfG_m)=\Ext_X^p(p^*\ZZ, \bfG_m)$,
 (see \cite[p. 23]{torsor} and \cite[Prop. 1.4.1]{colliot1987descente},
 respectively).

The existence of $\gamma$ follows from an easy diagram chase.
Since $u\cup\beta\in\im\Br_1X$,
\eqn{
u\cup\partial(\beta)=d(u\cup\beta)=0.
}
Then there exists $\gamma\in\Ext_k^2(S, \ol k[X]^\tm)$ such that
 $i^*(\gamma)=\partial(\beta)$.
Then  \eqref{cor_brextchase.it_chase} follows.
}

\rk{
We have some remarks on Corollary \ref{cor_brextchase}.
\enmt{[\upshape (1)]
\item One may use
  the variant Proposition \ref{prop_var_cone} to replace $\Pic\ol X$
  (resp. $\ol k[X]^\tm$) appearing in the diagrams  by
  $\KD'(X)$ (resp. $\ol k^\tm$).
 To be precise, take $f$ in Proposition \ref{prop_var_cone}  to be
  $\bfG_{m,k}\ra \tau_{\le1}\bfR p_*\bfG_{m,X}$. Then $\Delta=\KD'(X)$
  (see \cite{hs13descent}).
\item If moreover $M$ is finitely generated (hence so is $S$), then
 its Catier dual $\hat M$ is a group of multiplicative type.
 It can be shown that
\gan{
H^p(X, \hat T)=\Ext_X^p(p^*T, \bfG_m), \\
H^p(k, \hat T)=\Ext_k^p(T, \bfG_m)
}
 for $p\ge0$  and $T=M$ or $S$
  (c.f.  \cite[Prop. 1.4.1]{colliot1987descente}).
 Write $H_1^2(X, T) = \ker(H^2(X, T)\ra \Hom_k(T, \Br\ol X)$.
 Now  \eqref{eq_brextchase} becomes
\eqn{
\xymatrix@-2mm{
H^1(k, \ol k^\tm)\ar[r]\ar[d]_-{j^*} &\Pic X\ar[r]\ar[d]_-{p^*(j^*)}
  &\HH^0(k,\KD'(X))\ar[r]\ar[d]_-{j^*}
 &\Br k\ar[r]\ar[d]_-{j^*} &\Br_1X\ar[d]_-{p^*(j^*)} \\
H^1(k, \hat S) \ar[r]\ar[d]_-{i^*}
 &H^1(X, \hat S) \ar[r]^-\lambda\ar[d]_-{p^*(i^*)}
 &\Hom_{\bfD(k)}(S, \KD'(X))\ar[r]^-\partial\ar[d]_-{i^*}
 &H^2(k, \hat S)\ar[r]\ar[d]_-{i^*} &H^2_1(X, \hat S)
 \ar[d]_-{p^*(i^*)} \\
H^1(k, \hat M) \ar[r]\ar[d]_-{u\cup-}
 &H^1(X, \hat M) \ar[r]^-\lambda\ar[d]_-{p^*(u)\cup-}
 &\Hom_{\bfD(k)}(M, \KD'(X))\ar[r]^-\partial\ar[d]_-{u\cup-}
 &H^2(k, \hat M)\ar[r]\ar[d]_-{u\cup-} &H_1^2(X, \hat M)
 \ar[d]_-{p^*(u)\cup-} \\
\Br k\ar[r] &\Br_1X\ar[r]^-r &\HH^1(k, \KD'(X))\ar[r]^-d
 &H^3(k, \ol k^\tm)\ar[r] &H_{\le1}^3(X, \bfG_m)
}}
 where the second and third rows are \emph{fundamental exact sequences}
 for open varieties and the maps $\lambda$ are so-called \emph{extended
 type}. See \cite{hs13descent}.
}}


\bibliography{unibib}

\providecommand{\bysame}{\leavevmode\hbox to3em{\hrulefill}\thinspace}
\providecommand{\MR}{\relax\ifhmode\unskip\space\fi MR }
\providecommand{\MRhref}[2]{%
  \href{http://www.ams.org/mathscinet-getitem?mr=#1}{#2}
}
\providecommand{\href}[2]{#2}
\begin{thebibliography}{1}

\bibitem{colliot1987descente}
Jean-Louis Colliot-Th\'{e}l\`ene and Jean-Jacques Sansuc, \emph{La descente sur
  les vari\'{e}t\'{e}s rationnelles. {II}}, Duke Math. J. \textbf{54} (1987),
  no.~2, 375--492. \MR{899402}

\bibitem{gelfand2013methods}
S.~I. Gelfand and Y.~I. Manin, \emph{Methods of homological algebra}, Springer
  Science \& Business Media, 2013.

\bibitem{hs13descent}
David Harari and Alexei~N. Skorobogatov, \emph{Descent theory for open
  varieties}, Torsors, \'{e}tale homotopy and applications to rational points,
  London Math. Soc. Lecture Note Ser., vol. 405, Cambridge Univ. Press,
  Cambridge, 2013, pp.~250--279. \MR{3077172}

\bibitem{milne80etale}
J.~S. Milne, \emph{{\'E}tale cohomology}, Princeton Mathematical Series,
  vol.~33, Princeton Univ. Press, 1980.

\bibitem{skorobogatov99beyond}
Alexei~N. Skorobogatov, \emph{Beyond the {M}anin obstruction}, Invent. Math.
  \textbf{135} (1999), no.~2, 399--424. \MR{1666779}

\bibitem{torsor}
\bysame, \emph{Torsors and rational points}, vol. 144, Cambridge University
  Press, 2001.

\bibitem{skorobogatov2007eltob}
\bysame, \emph{On the elementary obstruction to the existence of rational
  points}, Mathematical Notes \textbf{81} (2007), no.~1-2, 97--107.

\bibitem{weibel1994homological}
C.~A. Weibel, \emph{An introduction to homological algebra}, Cambridge
  University Press, 1995.

\end{thebibliography}
\bibliographystyle{amsplain}
\end{document}